\documentclass{article}
\usepackage[intlimits]{amsmath}
\usepackage{amstext, amssymb}
\usepackage[mathscr]{euscript}
\usepackage{eufrak}
\usepackage[mathscr]{euscript}
\usepackage{cancel,enumerate}
\usepackage[affil-it]{authblk}
\usepackage{mathtools}

\title{Formula to evaluate 
$
\lim_{n\rightarrow \infty }\frac{1}{n}\sum_{i_{1},i_{2},\ldots, i_{k}=1}^{n}%
\lambda _{1}^{|i_{1}-i_{2}-s_{1}|}\lambda
_{2}^{|i_{2}-i_{3}-s_{2}|}\cdots\lambda _{k}^{|i_{k}-i_{1}-s_{k}|}
$
}
\date{\today}

\author[1]{Yuhao Liu}
\author[1]{Jan Vrbik}
\affil[1]{Department of Mathematics and Statistics\\ Brock University, Canada}

\begin{document}

\maketitle

\begin{abstract}
Computing moments of various parameter estimators related to an
autoregressive model of Statistics, one needs to evaluate several
expressions of the type mentioned in the title of this article. We proceed
to derive the corresponding formulas.
\end{abstract}

\section{Introduction}

The autoregressive model of Statistics generates a random sequence of observations by
\begin{equation}\label{ARk}
X_{i}=\alpha _{1}X_{i-1}+\alpha _{2}X_{i-2}+\cdots+\alpha
_{k}X_{i-k}+\varepsilon _{i}
\end{equation}
where $\varepsilon_{i}$ are independent, Normally distributed random variables with the mean of $0$ and the same standard deviation, and $k$ is a fixed integer, usually quite small (e.g. $k=1$ defines the so called Markov model). The sufficient and necessary condition for the resulting sequence to be asymptotically stationary is that all $k$ solutions of the characteristic polynomial
\begin{equation}
\lambda ^{k}=\alpha _{1}\lambda ^{k-1}+\alpha _{2}\lambda ^{k-2}+\cdots+\alpha_{k}  \label{CHP}
\end{equation}
are, in absolute value, smaller than $1$ (this is then assumed from now on).

The $j$th-order serial correlation coefficient $\rho _{j}$ (between $X_{i} $ and $X_{i+j}$) is then computed by
\begin{equation}\label{rho}
\rho _{j}=A_{1}\lambda _{1}^{|j|}+A_{2}\lambda _{2}^{|j|}+\cdots+A_{k}\lambda_{k}^{|j|}
\end{equation}%
where the $\lambda_{i}$'s are the $k$ roots of (\ref{CHP}), and the $A_{i}$ coefficients are themselves simple functions of these roots.

Computing the first few moments of various estimators (of the $\alpha _{i}$ parameters) boils down to computing moments of expressions of the
\begin{equation}\label{SX}
\sum_{i=1}^{n}X_{i}
\end{equation}
and 
\begin{equation}\label{DX}
\sum_{i=1}^{n-j}X_{i}X_{i+j}
\end{equation}
type, where $X_{1},$ $X_{2},\cdots X_{n}$ is a collection of $n$ consecutive observations (assuming that the process has already reached its stationary phase).

This in turn requires evaluating various summations (see \cite{vrbik}), of which the most difficult are
\begin{align}
\sum_{i_{1},i_{2}=1}^{\tilde{n}}\lambda _{1}^{|i_{1}-i_{2}+s_{1}|}\lambda_{2}^{|i_{2}-i_{1}+s_{2}|}  \label{S2} \\
\sum_{i_{1},i_{2},i_{3}=1}^{\tilde{n}}\lambda_{1}^{|i_{1}-i_{2}+s_{1}|}\lambda _{2}^{|i_{2}-i_{3}+s_{2}|}\lambda_{3}^{|i_{3}-i_{1}+s_{3}|}  \label{S3}
\end{align}
and
\begin{align}
\sum_{i_{1},i_{2},i_{3},i_{4}=1}^{\tilde{n}}\lambda
_{1}^{|i_{1}-i_{2}+s_{1}|}\lambda _{2}^{|i_{2}-i_{3}+s_{2}|}\lambda_{3}^{|i_{3}-i_{4}+s_{3}|}\lambda _{4}^{|i_{4}-i_{1}+s_{4}|}  \label{S4}
\end{align}
where $\lambda_{1},$ $\lambda_{2},$ $\lambda_{3},$ and $\lambda _{4}$ are any of the $\lambda _{i}$ roots (some may appear in duplicate), $s_{1}$, $s_{2}$, $s_{3}$, and $s_{4}$ are small integers, and $\tilde{n}$ indicates that each of the upper limits equals to $n$, perhaps adjusted in the manner of (\ref{DX}).

It is possible (but rather messy --- the result depends on the values of $s_{1} $, $s_{2}$ and $\tilde{n}$ --- see \cite{LIU}) to \emph{exactly} evaluate (\ref{S2}) and realize that the answer will always (this goes for the other two summations as well) consist of three parts:

\begin{enumerate}
\item terms proportional to $\lambda _{i}^{n},$ which tend to zero (as $n$ increases) `exponentially',
\item terms which stay constant as $n$ increases,
\item terms proportional to $n.$
\end{enumerate}

Luckily, to build the approximation which is usually deemed sufficient (see  \cite{vrbik}), we need to find only the $n$ proportional terms. These can be extracted by dividing the relevant summation by $n$ and taking the $n \to \infty$ limit. Incidentally, this results in the following (and most welcomed) simplification: the corresponding answer will be the  \emph{same} regardless of the $\tilde{n}$ adjustments (thus, we may as well use $n$ instead), and will similarly \emph{not} depend on the individual $s_{i}$'s, but only on the \emph{absolute} value of their \emph{sum}. The proof of this statement is omitted.

\section{Evaluating the limits}

Starting with (\ref{S2}), we obtain
\begin{align*}
\MoveEqLeft F_{2}(\lambda _{1},\lambda _{2};S)\\
&\equiv \lim_{n\rightarrow \infty }\frac{1}{n}\sum_{i_{1},i_{2}=1}^{\tilde{n}}\lambda
_{1}^{|i_{1}-i_{2}+s_{1}|}\lambda _{2}^{|i_{2}-i_{1}+s_{2}|} \\
&=\lim_{n\rightarrow \infty }\frac{1}{n}\sum_{i_{1},i_{2}=1}^{n}\lambda _{1}^{|i_{1}-i_{2}|}\lambda_{2}^{|i_{2}-i_{1}+S|}=\lim_{n\rightarrow \infty }\frac{1}{n}\sum_{j=-n}^{n}\sum_{i_{1}=\max (1,1-j)}^{\min (n,n-j)}\lambda_{1}^{|j|}\lambda _{2}^{|S-j|} \\
&=\lim_{n\rightarrow \infty }\sum_{j=-n}^{n}\frac{n-|j|}{n}\lambda _{1}^{|j|}\lambda _{2}^{|S-j|}=\sum_{j=-\infty }^{\infty}\lambda _{1}^{|j|}\lambda _{2}^{|S-j|} \\
&=\sum_{j=-\infty }^{0}\lambda _{1}^{-j}\lambda_{2}^{S-j}+\sum_{j=1}^{S}\lambda _{1}^{j}\lambda
_{2}^{S-j}+\sum_{j=S+1}^{\infty }\lambda _{1}^{j}\lambda _{2}^{j-S} \\
&=\frac{\lambda _{1}^{S+1}(1-\lambda _{2}^{2})}{(\lambda_{1}-\lambda _{2})(1-\lambda _{1}\lambda _{2})}+\frac{\lambda_{2}^{S+1}(1-\lambda _{1}^{2})}{(\lambda _{2}-\lambda _{1})(1-\lambda_{2}\lambda _{1})}
\end{align*}
where $S\equiv |s_{1}+s_{2}|$. Following the usual convention, an empty summation (such as $\sum_{j=1}^{0}$) has a zero value.

Note that the answer can be written in the following form:%
\begin{equation}
\sum_{i=1}^{\ell }\lambda _{i}^{S+\ell -1}\prod_{j=1,j\neq i}^{\ell }\frac{%
1-\lambda _{j}^{2}}{(\lambda _{i}-\lambda _{j})(1-\lambda _{i}\lambda _{j})}
\label{MAIN}
\end{equation}%
with $\ell =2.$ Also note that, when $\lambda _{2}=\lambda _{1},$ the value
of $F_{2}(\lambda _{1},\lambda _{1};S)$ can be easily obtained by%
\begin{equation*}
\lim_{\lambda _{2}\rightarrow \lambda _{1}}F_{2}(\lambda _{1},\lambda
_{2};S)=\frac{\lambda _{1}^{S}\left( 1+S\overset{}{+}(1-S)\lambda
_{1}^{2}\right) }{1-\lambda _{1}^{2}}
\end{equation*}

\subsection{The case of 3 $\protect\lambda $'s}

Moving on to (\ref{S3}), we now get%
\begin{align*}
\MoveEqLeft F_{3}(\lambda _{1},\lambda _{2},\lambda _{3};S)\equiv \lim_{n\rightarrow
\infty }\frac{1}{n}\sum_{i_{1},i_{2},i_{3}=1}^{\tilde{n}}\lambda
_{1}^{|i_{1}-i_{2}+s_{1}|}\lambda _{2}^{|i_{2}-i_{3}+s_{2}|}\lambda
_{3}^{|i_{3}-i_{1}+s_{3}|} \\
&=\lim_{n\rightarrow \infty }\frac{1}{n}%
\sum_{i_{1},i_{2},i_{3}=1}^{n}\lambda _{1}^{|i_{1}-i_{2}|}\lambda
_{2}^{|i_{2}-i_{3}|}\lambda _{3}^{|i_{3}-i_{1}+S|} \\
&=\lim_{n\rightarrow \infty }\frac{1}{n}\sum_{j_{1}=-n}^{n}%
\sum_{j_{2}=\max (-n,-n-j_{1})}^{\min (n,n-j_{1})}\sum_{i_{1}=\max
(1,1-j_{2},1-j_{1}-j_{2})}^{\min (n,n-j_{2},n-j_{1}-j_{2})}\lambda
_{1}^{|j_{1}|}\lambda _{2}^{|j_{2}|}\lambda _{3}^{|S-j_{1}-j_{2}|} \\
&=\lim_{n\rightarrow \infty
}\sum_{j_{1}=-n}^{n}\sum_{j_{2}=\max (-n,-n-j_{1})}^{\min (n,n-j_{1})}\tfrac{%
\min (n,n-j_{2},n-j_{1}-j_{2})-\max (1,1-j_{2},1-j_{1}-j_{2})+1}{n}\cdot
\lambda _{1}^{|j_{1}|}\lambda _{2}^{|j_{2}|}\lambda _{3}^{|S-j_{1}-j_{2}|} \\
&=\sum_{j_{1}=-\infty }^{\infty }\sum_{j_{2}=-\infty }^{\infty
}\lambda _{1}^{|j_{1}|}\lambda _{2}^{|j_{2}|}\lambda _{3}^{|S-j_{1}-j_{2}|}
\end{align*}%
where $S\equiv |s_{1}+s_{2}+s_{3}|.$

This time, evaluating the last summation is slightly more difficult; we will
do it quadrant by quadrant.

For the first quadrant (\emph{including} the adjacent half-axes and the
origin), we get (visualize the quadrant, cut by the $S=j_{1}+j_{2}$ line):%
\begin{align*}
\MoveEqLeft 
\sum_{j_{1}=0}^{S}\sum_{j_{2}=0}^{S-j_{1}}\lambda _{1}^{j_{1}}\lambda
_{2}^{j_{2}}\lambda
_{3}^{S-j_{1}-j_{2}}+\sum_{j_{1}=0}^{S}\sum_{j_{2}=S-j_{1}+1}^{\infty
}\lambda _{1}^{j_{1}}\lambda _{2}^{j_{2}}\lambda
_{3}^{j_{1}+j_{2}-S}+\sum_{j_{1}=S+1}^{\infty }\sum_{j_{2}=0}^{\infty
}\lambda _{1}^{j_{1}}\lambda _{2}^{j_{2}}\lambda _{3}^{j_{1}+j_{2}-S} 
\\
&= \frac{\lambda _{1}^{S+2}(1-\lambda _{3}^{2})}{(\lambda _{1}-\lambda
_{2})(\lambda _{1}-\lambda _{3})(1-\lambda _{1}\lambda _{3})}+\frac{\lambda
_{2}^{S+2}(1-\lambda _{3}^{2})}{(\lambda _{2}-\lambda _{1})(\lambda
_{2}-\lambda _{3})(1-\lambda _{2}\lambda _{3})}+\frac{\lambda _{3}^{S+2}}{%
(\lambda _{3}-\lambda _{1})(\lambda _{3}-\lambda _{2})}
\end{align*}%
For the second quadrant (again, \emph{including} the corresponding
boundaries - the resulting duplication with the first quadrant will be
removed later), the same kind of approach yields%
\begin{align*}
\MoveEqLeft \sum_{j_{1}=-\infty }^{0}\sum_{j_{2}=0}^{S-j_{1}}\lambda
_{1}^{-j_{1}}\lambda _{2}^{j_{2}}\lambda
_{3}^{S-j_{1}-j_{2}}+\sum_{j_{1}=-\infty }^{0}\sum_{j_{2}=S-j_{1}+1}^{\infty
}\lambda _{1}^{-j_{1}}\lambda _{2}^{j_{2}}\lambda _{3}^{j_{1}+j_{2}-S} \\
&=\frac{\lambda _{2}^{S+1}(1-\lambda _{3}^{2})}{(\lambda _{2}-\lambda
_{3})(1-\lambda _{2}\lambda _{1})(1-\lambda _{2}\lambda _{3})}+\frac{\lambda
_{3}^{S+1}}{(\lambda _{3}-\lambda _{2})(1-\lambda _{3}\lambda _{1})}
\end{align*}%
The fourth quadrant clearly results in the same answer, with $\lambda _{1}$
and $\lambda _{2}$ interchanged, namely%
\begin{equation*}
\frac{\lambda _{1}^{S+1}(1-\lambda _{3}^{2})}{(\lambda _{1}-\lambda
_{3})(1-\lambda _{1}\lambda _{2})(1-\lambda _{1}\lambda _{3})}+\frac{\lambda
_{3}^{S+1}}{(\lambda _{3}-\lambda _{1})(1-\lambda _{3}\lambda _{2})}
\end{equation*}%
Finally, the third quadrant (including its boundaries) contributes%
\begin{align*}
\sum_{j_{1}=-\infty }^{0}\sum_{j_{2}=-\infty }^{0}\lambda
_{1}^{-j_{1}}\lambda _{2}^{-j_{2}}\lambda _{3}^{S-j_{1}-j_{2}} 
&=\frac{\lambda _{3}^{S}}{(1-\lambda _{3}\lambda _{1})(1-\lambda
_{3}\lambda _{2})}.
\end{align*}

Adding the four results does \emph{not} yield the desired answer, since the
contribution of each of the two axes has been included \emph{twice}, and
that of the origin altogether \emph{four} times. This can be easily
corrected by subtracting $F_{2}(\lambda _{2},\lambda _{3};S)$ which removes
the extra contribution of the $j_{1}=0$ axis, and $F_{2}(\lambda
_{1},\lambda _{3};S)$ which does the same thing with the $j_{2}=0$ terms.
This leaves us with the origin ($j_{1}=j_{2}=0$) which, at this point, is
still contributing \emph{double} its value (two contributions have been
removed with the two axes); subtracting $\lambda _{3}^{S}$ fixes that as
well.

The final answer thus becomes
\begin{align*}
F_{3}(\lambda _{1},\lambda _{2},\lambda _{3};S) 
=~&
\frac{\lambda
_{1}^{S+2}(1-\lambda _{2}^{2})(1-\lambda _{3}^{2})}{(\lambda _{1}-\lambda
_{2})(\lambda _{1}-\lambda _{3})(1-\lambda _{1}\lambda _{2})(1-\lambda
_{1}\lambda _{3})} \\
&+  \frac{\lambda _{2}^{S+2}(1-\lambda _{1}^{2})(1-\lambda _{3}^{2})}{(\lambda
_{2}-\lambda _{1})(\lambda _{2}-\lambda _{3})(1-\lambda _{2}\lambda
_{1})(1-\lambda _{2}\lambda _{3})} \\
&+\frac{\lambda _{3}^{S+2}(1-\lambda _{1}^{2})(1-\lambda _{2}^{2})}{(\lambda
_{3}-\lambda _{1})(\lambda _{3}-\lambda _{2})(1-\lambda _{3}\lambda
_{1})(1-\lambda _{3}\lambda _{2})}
\end{align*}%
Note that this has the form of (\ref{MAIN}) with $\ell =3.$

When any of the three $\lambda $`s are identical, the answer can be found as
the corresponding limit of the previous expression. Thus, for example%
\begin{equation*}
F_{3}(\lambda ,\lambda ,\lambda ;S)=\frac{2+3S+S^{2}+2(4-S^{2})\lambda
^{2}+(2-3S+S^{2})\lambda ^{4}}{2(1-\lambda ^{2})^{2}}\cdot \lambda ^{S}
\end{equation*}%
etc.

\subsection{The case of 4 $\protect\lambda $'s}

The main challenge is to evaluate the last limit, namely%
\begin{align*}
\MoveEqLeft F_{4}(\lambda _{1},\lambda _{2},\lambda _{3},\lambda _{4};S) \\
&\equiv\lim_{n\rightarrow \infty }\frac{1}{n}\sum_{i_{1},i_{2},i_{3},i_{4}=1}^{\tilde{n}}\lambda _{1}^{|i_{1}-i_{2}+s_{1}|}\lambda
_{2}^{|i_{2}-i_{3}+s_{2}|}\lambda _{3}^{|i_{3}-i_{4}+s_{3}|}\lambda
_{4}^{|i_{4}-i_{1}+s_{4}|} \\
&=\sum_{j_{1}=-\infty }^{\infty }\sum_{j_{2}=-\infty }^{\infty
}\sum_{j_{3}=-\infty }^{\infty }\lambda _{1}^{|j_{1}|}\lambda
_{2}^{|j_{2}|}\lambda _{3}^{|j_{3}|}\lambda _{4}^{|S-j_{1}-j_{2}-j_{3}|}
\end{align*}%
We proceed octant by octant; the octants will be identified by the signs of
the $j_{1},$ $j_{2}$ and $j_{3}$ indices, respectively.

For the first octant denoted $O_{+++}$ (including the adjacent portions of
coordinate planes, axes and the origin), we get
\begin{align*}
\MoveEqLeft
\sum_{j_{1}=0}^{S}\sum_{j_{2}=0}^{S-j_{1}}\sum_{j_{3}=0}^{S-j_{1}-j_{2}}
\lambda _{1}^{j_{1}}\lambda _{2}^{j_{2}}\lambda _{3}^{j_{3}}\lambda
_{4}^{S-j_{1}-j_{2}-j_{3}} \\
&+\sum_{j_{1}=0}^{S}\sum_{j_{2}=0}^{S-j_{1}}\sum_{j_{3}=S-j_{1}-j_{2}+1}^{\infty }\lambda _{1}^{j_{1}}\lambda_{2}^{j_{2}}\lambda _{3}^{j_{3}}\lambda _{4}^{j_{1}+j_{2}+j_{3}-S} \\
&+\sum_{j_{1}=0}^{S}\sum_{j_{2}=S-j_{1}+1}^{\infty }\sum_{j_{3}=0}^{\infty}\lambda _{1}^{j_{1}}\lambda _{2}^{j_{2}}\lambda _{3}^{j_{3}}\lambda_{4}^{j_{1}+j_{2}+j_{3}-S}\\
&+\sum_{j_{1}=S+1}^{\infty }\sum_{j_{2}=0}^{\infty}\sum_{j_{3}=0}^{\infty }\lambda _{1}^{j_{1}}\lambda _{2}^{j_{2}}\lambda_{3}^{j_{3}}\lambda _{4}^{j_{1}+j_{2}+j_{3}-S} \\
=&~
\frac{\lambda _{1}^{S+3}(1-\lambda _{4}^{2})}{(\lambda _{1}-\lambda
_{2})(\lambda _{1}-\lambda _{3})(\lambda _{1}-\lambda _{4})(1-\lambda
_{1}\lambda _{4})} \\
&+\frac{\lambda _{2}^{S+3}(1-\lambda _{4}^{2})}{(\lambda_{2}-\lambda _{1})(\lambda _{2}-\lambda _{3})(\lambda _{2}-\lambda_{4})(1-\lambda _{2}\lambda _{4})} \\
&+\frac{\lambda _{3}^{S+3}(1-\lambda _{4}^{2})}{(\lambda _{3}-\lambda
_{2})(\lambda _{3}-\lambda _{2})(\lambda _{3}-\lambda _{4})(1-\lambda
_{3}\lambda _{4})} \\
&+\frac{\lambda _{4}^{S+3}}{(\lambda _{4}-\lambda_{1})(\lambda _{4}-\lambda _{2})(\lambda _{4}-\lambda _{3})}
\end{align*}%
To understand why it was necessary to break the summation into four parts,
it helps to visualize the first octant, cut by the $S=j_{1}+j_{2}+j_{3}$
plane,thus:
$$
\includegraphics[width=.6\textwidth]{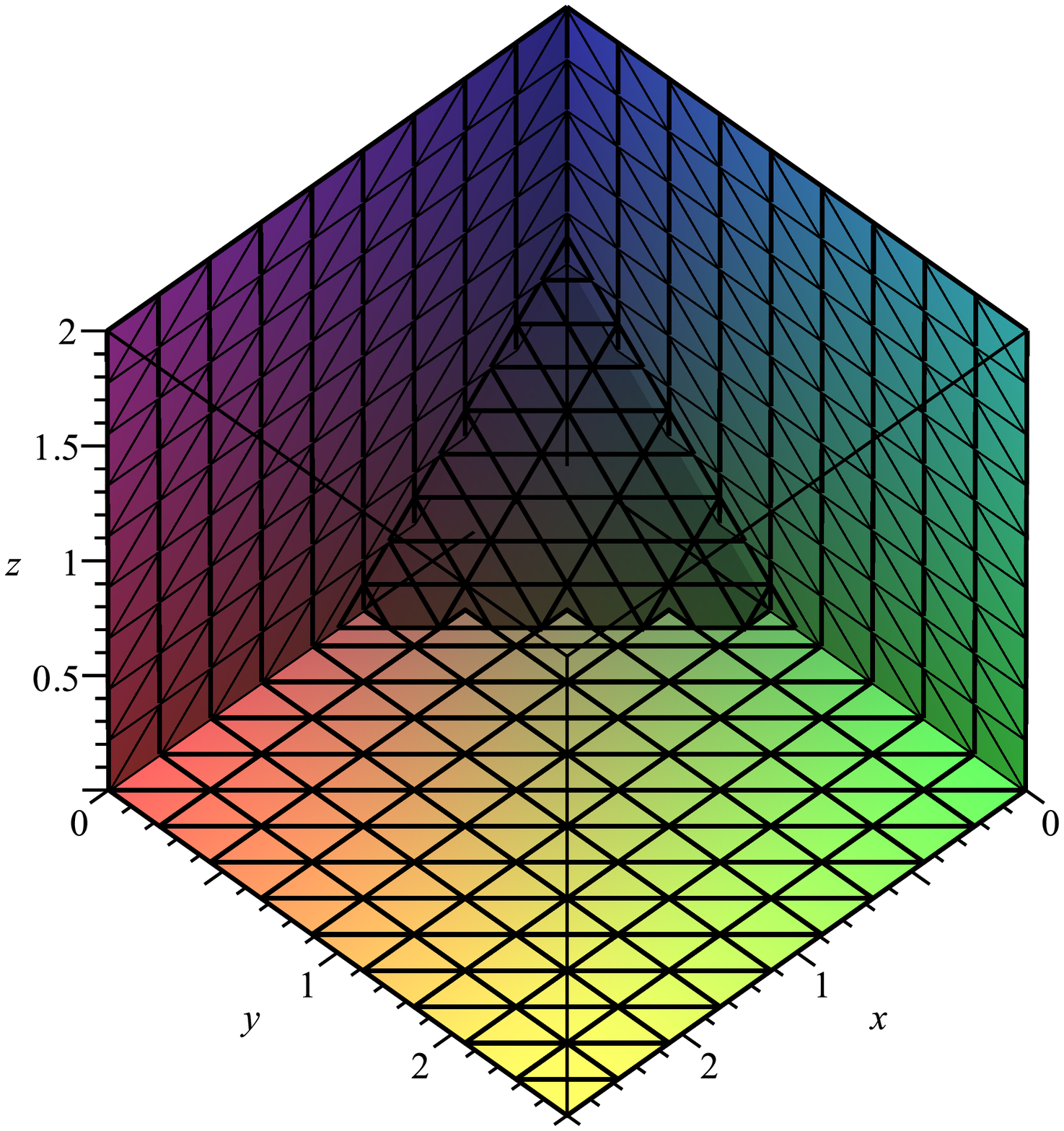}
$$
Our brain can interpret this image in two different ways; please make an
effort to see the triangle as the most \emph{distant} part of the picture.

As the next octant we take $O_{++-}$ (with all its boundaries), contributing
\begin{align*}
&
\sum_{j_{1}=0}^{S}\sum_{j_{2}=0}^{S-j_{1}}\sum_{j_{3}=-\infty }^{0}\lambda_{1}^{j_{1}}\lambda _{2}^{j_{2}}\lambda _{3}^{-j_{3}}\lambda
_{4}^{S-j_{1}-j_{2}-j_{3}}\\
&+\sum_{j_{1}=0}^{\infty}\sum_{j_{2}=S-j_{1}+1}^{\infty }\sum_{j_{3}=-\infty}^{S-j_{1}-j_{2}}\lambda _{1}^{j_{1}}\lambda _{2}^{j_{2}}\lambda_{3}^{-j_{3}}\lambda _{4}^{S-j_{1}-j_{2}-j_{3}} \\
&+\sum_{j_{1}=0}^{\infty }\sum_{j_{2}=S-j_{1}+1}^{\infty}\sum_{j_{3}=S-j_{1}-j_{2}+1}^{0}\lambda _{1}^{j_{1}}\lambda_{2}^{j_{2}}\lambda _{3}^{-j_{3}}\lambda_{4}^{j_{1}+j_{2}+j_{3}-S}\\
&+\sum_{j_{1}=S+1}^{\infty }\sum_{j_{2}=0}^{\infty}\sum_{j_{3}=-\infty }^{S-j_{1}-j_{2}}\lambda _{1}^{j_{1}}\lambda_{2}^{j_{2}}\lambda _{3}^{-j_{3}}\lambda _{4}^{S-j_{1}-j_{2}-j_{3}} \\
&+\sum_{j_{1}=S+1}^{\infty }\sum_{j_{2}=0}^{\infty}\sum_{j_{3}=S-j_{1}-j_{2}+1}^{0}\lambda _{1}^{j_{1}}\lambda_{2}^{j_{2}}\lambda _{3}^{-j_{3}}\lambda _{4}^{j_{1}+j_{2}+j_{3}-S}\\
&=
\frac{\lambda _{1}^{S+2}(1-\lambda _{4}^{2})}{(\lambda _{1}-\lambda_{2})(\lambda _{1}-\lambda _{4})(1-\lambda _{1}\lambda _{3})(1-\lambda_{1}\lambda _{4})} \\
&~~~+\frac{\lambda _{2}^{S+2}(1-\lambda _{4}^{2})}{(\lambda_{2}-\lambda _{1})(\lambda _{2}-\lambda _{4})(1-\lambda _{2}\lambda_{3})(1-\lambda _{2}\lambda _{4})} \\
&~~~+\frac{\lambda _{4}^{S+2}}{(\lambda _{4}-\lambda _{1})(\lambda_{4}-\lambda _{2})(1-\lambda _{4}\lambda _{3})}
\end{align*}
Again, visualizing the situation may help (the corner being the most \emph{distant} part of the picture):
$$
\includegraphics[width=.6\textwidth]{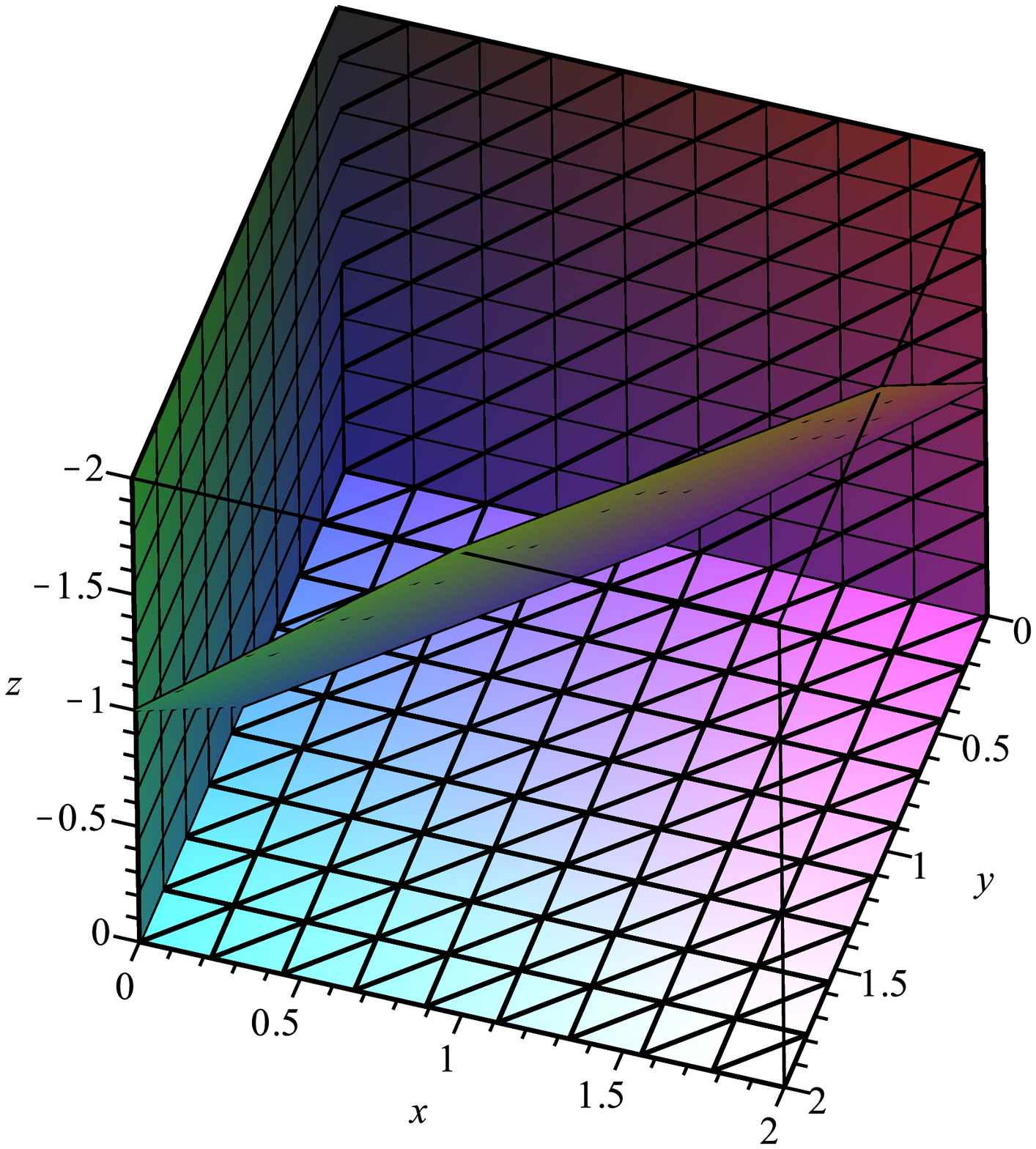}
$$

The $O_{-++}$ and $O_{+-+}$ octants contribute the same expression each,
after the $\lambda _{3}\leftrightarrow \lambda _{1}$ and $\lambda
_{3}\leftrightarrow \lambda _{2}$ interchange, respectively.

For $O_{--+}$ (including boundaries) we get
\begin{align*}
\MoveEqLeft 
\sum_{j_{1}=-\infty }^{0}\sum_{j_{2}=-\infty}^{0}\sum_{j_{3}=0}^{S-j_{1}-j_{2}}\lambda _{1}^{-j_{1}}\lambda_{2}^{-j_{2}}\lambda _{3}^{j_{3}}\lambda_{4}^{S-j_{1}-j_{2}-j_{3}}\\
&+\sum_{j_{1}=-\infty }^{0}\sum_{j_{2}=-\infty}^{0}\sum_{j_{3}=S-j_{1}-j_{2}+1}^{\infty }\lambda _{1}^{-j_{1}}\lambda_{2}^{-j_{2}}\lambda _{3}^{j_{3}}\lambda _{4}^{j_{1}+j_{2}+j_{3}-S} \\[.4em]
=~& \frac{\lambda _{3}^{S+1}(1-\lambda _{4}^{2})}{(\lambda _{3}-\lambda_{4})(1-\lambda _{3}\lambda _{1})(1-\lambda _{3}\lambda _{2})(1-\lambda_{3}\lambda_{4})}\\
&+\frac{\lambda _{4}^{S+1}}{(\lambda _{4}-\lambda_{3})(1-\lambda _{4}\lambda _{1})(1-\lambda _{4}\lambda _{2})}
\end{align*}%
because this is how it looks like (again, the corner to be seen as most
distant)%
$$
\includegraphics[width=.6\textwidth]{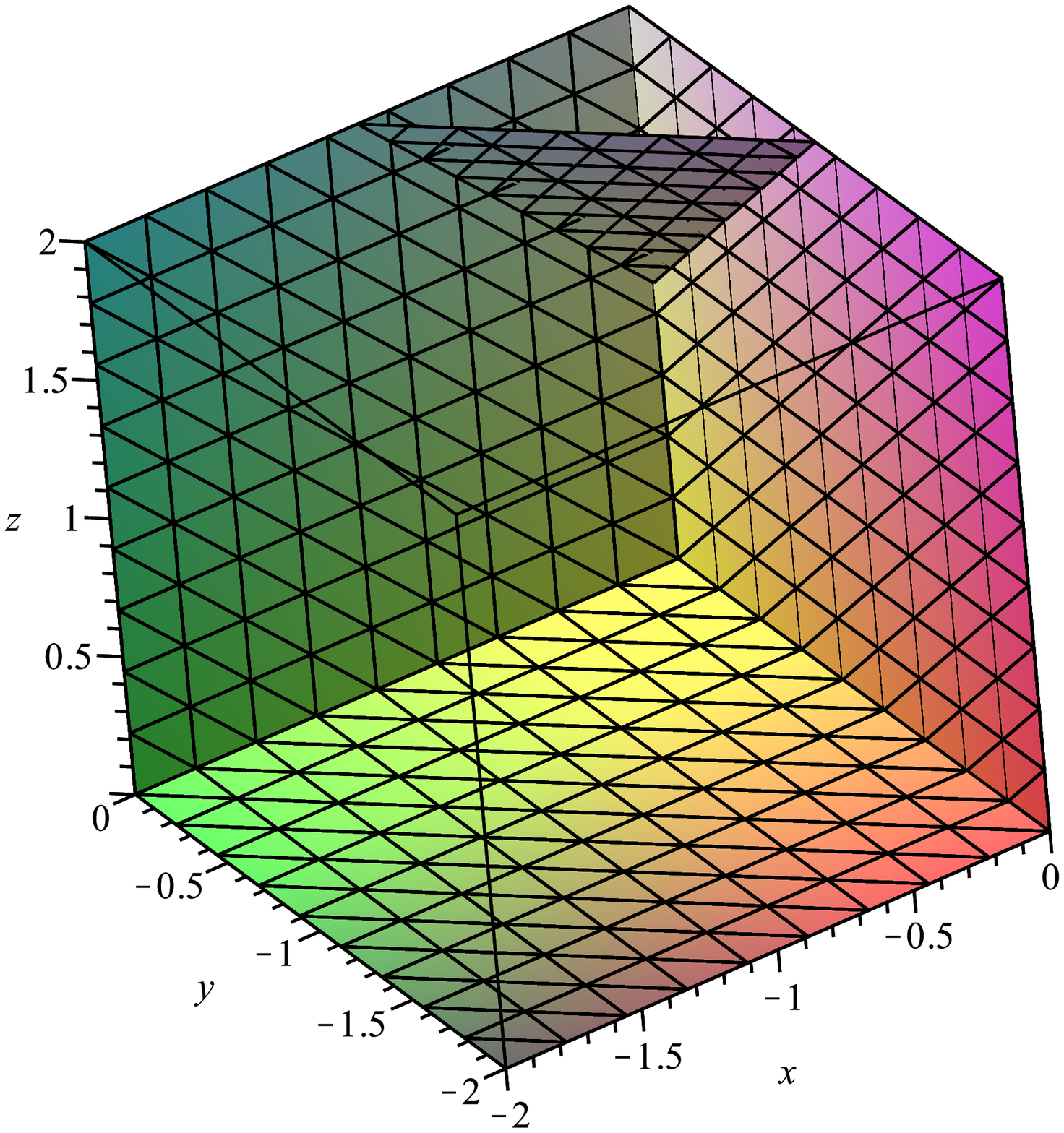}
$$
and similarly for $O_{+--}$ and $O_{-+-,}$ after the $\lambda
_{3}\leftrightarrow \lambda _{1}$ and $\lambda _{3}\leftrightarrow \lambda
_{2}$ interchange, respectively.

Finally, $O_{\_\_\_}$ with its boundaries contributes%
\begin{align*}
\MoveEqLeft 
\sum_{j_{1}=-\infty }^{0}\sum_{j_{2}=-\infty }^{0}\sum_{j_{3}=-\infty
}^{0}\lambda _{1}^{-j_{1}}\lambda _{2}^{-j_{2}}\lambda _{3}^{-j_{3}}\lambda
_{4}^{S-j_{1}-j_{2}-j_{3}} \\
&=\frac{\lambda _{4}^{S}}{(1-\lambda _{4}\lambda _{1})(1-\lambda_{4}\lambda _{2})(1-\lambda _{4}\lambda _{3})}
\end{align*}

Adding the eight results and subtracting $F_{3}(\lambda _{2},\lambda
_{3},\lambda _{4};S)+F_{3}(\lambda _{1},\lambda _{3},\lambda
_{4};S)+F_{3}(\lambda _{1},\lambda _{2},\lambda _{4};S)$ to remove the
duplicate contribution of the three coordinate planes; further subtracting $%
F_{2}(\lambda _{1},\lambda _{4};S)+F_{2}(\lambda _{2},\lambda
_{4};S)+F_{2}(\lambda _{3},\lambda _{4};S)$ to remove the originally
quadruple (now duplicate) contribution of the three axes; and finally
subtracting $\lambda _{4}^{S}$ to remove the remaining, originally eightfold
(now duplicate) contribution of the origin, yields the final formula for $%
F_{4}(\lambda _{1},\lambda _{2},\lambda _{3},\lambda _{4};S).$ Not
surprisingly, it turns out to be equal to (\ref{MAIN}) with $\ell =4.$

\subsection{Further challenge}

At this point, it is fairly obvious that $F_{5}(\lambda _{1},\lambda
_{2},\lambda _{3},\lambda _{4},\lambda _{5};S)$ will be given by (\ref{MAIN}%
) with $\ell =5,$ etc. To prove this by the technique of this article
becomes increasingly more difficult (impossible in general, since $\ell $
can have any integer value). One clearly needs to proceed by induction -
would anyone want to try?

\end{document}